\documentclass{amsart} 
 \newcommand\Pn[1]{{\bf P}^{#1}} 
 \newcommand\Pnd[1]{{\check{\bf P}}^{#1}}
 \newtheorem{theorem}{{\rm T\sc heorem}}[section] 
 \newtheorem{lemma}[theorem]{{\rm L\sc emma}} 
 \newtheorem{corollary}[theorem]{{\rm C\sc orollary}} 
 \newtheorem{proposition}[theorem]{{\rm P\sc roposition}}

 \newtheorem{problem}[theorem]{Problem} 
  
 \newtheorem{re}[theorem]{{\rm R\sc emark}}

 \newcommand\cF{\mathcal F} 
 \newcommand\C{{\bf C}} 
 \newcommand\Gr{{\bf G}}
 
 \newcommand\Hom{\mathop{\rm Hom}\nolimits} 
  
 \newcommand\Pfaff{\mathop{\rm Pfaff}\nolimits} 
 \newcommand\rto{--->}

\begin{document}
 \title {${\bf K3}$ surfaces of genus 8 and varieties of sums of powers of cubic 
 fourfolds}
 \author{Atanas Iliev and Kristian Ranestad}
 
 \begin{abstract} 
The main outcome of this paper is that the variety VSP(F,10) of presentations 
of a general cubic form F in 6 variables as a sum of 10 cubes is a smooth 
symplectic 4-fold obtained a deformation of the Hilbert square of a K3 surface 
of genus 8. After publishing it in Trans. Am. Math. Soc. 353, No.4, 
1455-1468 (2001), it was noted to us by Eyal Markman that in Theorem 3.17 
we conclude without proof that VSP(F,10) should be the 4-fold of lines on 
another cubic 4-fold. 
We correct this in the e-print "Addendum to K3 surfaces of genus 8 and 
varieties of sums of powers of cubic fourfolds" (math.AG/0611533), where we establish
that the general VSP(F,10) is in fact a new symplectic 4-fold different from 
the family of lines on a cubic 4-fold. 
 \end{abstract}
 \maketitle 
 {\footnotetext{Mathematics Subject Classification 14J70 (Primary), 14M15, 
 14N99 (Secondary)}}
 \pagestyle{myheadings}\markboth{\textsc{Atanas Iliev and Kristian 
 Ranestad}}{\textsc{K3 surfaces of genus 8 and varieties of sums of powers 
 of cubic 
 4-folds}}
 \section{Pfaffian and apolar cubic 4-folds associated to $K3$ surfaces 
 of genus 8}{\label{s1}}
 \subsection{\label{1.1}} Let $V$ be a 6-dimensional vector space over $\C$. 
 Fix a basis $e_0,\ldots ,e_5$ for $V$, then $e_i\wedge e_j$ for 
 $0\leq i<j\leq 5$ form a basis for the Pl\"ucker space $\wedge^2V$ of 
 2-dimensional subspaces in 
 $V$ or lines in $\Pn 5={\bf P}(V)$. We associate to a $2$-vector 
 $g=\sum_{i<j}a_{ij}e_{i}\wedge e_{j}\in \wedge^2V$ a skewsymmetric 
 matrix $M(g)=(a_{ij})$, with $a_{ji}=-a_{ij}$. With Pl\"ucker coordinates 
 $x_{ij}$, 
 the embedding of the Grassmannian $G=\Gr (2,V)$ in $\Pn 
 {14}={\bf P}(\Lambda^2V)$ is then precisely the locus of rank 2 
 skewsymmetric $6\times 6$ matrices 
 $$M=\left(\begin{array}{cccccc} 0&x_{01}&x_{02}&x_{03}&x_{04}&x_{05}\\ 
 -x_{01}&0&x_{12}&x_{13}&x_{14}&x_{15}\\ 
 -x_{02}&-x_{12}&0&x_{23}&x_{24}&x_{25}\\ 
 -x_{03}&-x_{13}&-x_{23}&0&x_{34}&x_{35}\\ 
 -x_{04}&-x_{14}&-x_{24}&-x_{34}&0&x_{45}\\ 
 -x_{05}&-x_{15}&-x_{25}&-x_{35}&-x_{45}&0 
 \end{array}\right) 
 $$
 Since the 
 sum of two rank 2 matrices has rank at most 4, and any rank 4 skewsymmetric 
 matrix is the 
 sum of two rank 2 skewsymmetric matrices, the secant variety of $G$ is the 
 cubic 
 hypersurface $K$ defined by the $6\times 6$ Pfaffian $m$ of the matrix $M$. 
 The dual 
 variety of $G$ in $\Pnd {14}={\bf P}(\wedge^2V^{*})$ is a 
 cubic hypersurface $K^*\cong K$ (cf. \cite{Zak}). $K^*$ is the secant 
 variety of $G^*=\Gr (V,2)$ the Grassmannian of rank $2$ quotient spaces of 
 $V$, and of course 
 $G^*\cong G$.
 \subsection{\label{1.2}} 
 A general $K3$ surface $S$ with 
 Picard group generated by a linebundle $H$ of degree $H^2=14$ is embedded 
 via $|H|$ into the 
 Grassmannian $G^*=\Gr (V,2)$. In 
 fact $S$ is the intersection of $G^*$ with a linear space $L_S$ of 
 dimension $8$ in $\Pnd 
 {14}={\bf P}(\wedge^2V^*)$ (cf. \cite{Muk}). The dual space 
 $P_S=L_S^{\bot}\subset \Pn 
 {14}$ is $5$-dimensional, so $P_S$ intersects the dual variety $K$ of $G^*$ 
 in a Pfaffian cubic 
 4-fold which we denote by $F'(S)$ (or simply $F'$).
 
 \subsection{\label{1.3}} The Pl\"ucker embedding of the Grassmannian $G=\Gr 
 (2,V)$ in 
 $\Pn {14}$ is arithmetically Gorenstein. The homogeneous coordinate ring 
 $A_{G}$ has 
 syzygies, easily computed with \cite{MAC},
 
 \[ \begin{array}{ccccccc} 
 1 & - & - & - & - & -& - \\ 
 - & 15 & 35 & 21& - & -& - \\ 
 - & - & - & 21& 35 & 15& - \\ 
 - & - & - & - & - & -& 1 
 \end{array} 
 \] 
 in MACAULAY notation. The Grassmannian variety has dimension $8$, so 
 $P_S=L_{S}^{\bot}$, defined by the linear forms $h_0,\ldots, h_8$ in $\Pn 
 {14}$, does not 
 intersect $G$, and the quotient $A=A_{G}/(h_0,\dots,h_8)$ is an Artinian 
 Gorenstein ring. Its 
 Hilbert function is $(1,6,6,1)$ with socledegree 3, so $A$ is the apolar 
 Artinian Gorenstein ring 
 $A^{F}$ for some cubic hypersurface $F=F(S)\subset \check P_S$. Thus the 
 dual socle 
 generator of $A$ is a cubic form $f$, defined up to scalar by $A$, 
 with $F=Z(f)$. 
 We call $F=F(S)$ the apolar cubic 4-fold of $S$.
 
 \setcounter{theorem}{3} 
 \begin{lemma}{\label{1.4}} There is a 19-dimensional family of cubic 
 $4$-folds $F$ whose apolar 
 Artinian Gorenstein ring is a quotient of $A_{G}$. 
 \end{lemma} 
 \proof Macaulay showed that there 
 is a $1:1$ correspondence between hypersurfaces of degree $d$ and graded 
 Artinian Gorenstein 
 rings generated in degree 1 with socledegree $d$ \cite{Mac} (cf. 
 \ref{3.3} 
 below). Now, an 
 isomorphism between such rings is of course induced by a linear 
 transformation on the 
 generators. In our setting any such linear transformation is again induced 
 by an automorphism 
 of $G^*$ and correspondingly of $G$. The isomorphism classes of general 
 $K3$ surfaces of 
 genus 8 correspond precisely to orbits of $8$-dimensional subspaces $L$ 
 (cf. \cite{Muk}). There is a $19$-dimensional family of 
 $K3$ surfaces of genus $8$, so the lemma follows.\qed\vspace{.2 in}

 \begin{re}{\label{1.5}} The apolar cubic $F=F(S)$ is in 
 general not a Pfaffian 
 cubic.\end{re} 
 In fact the Pfaffian cubics 
 also form a $19$-dimensional family of cubic $4$-folds, and the above 
 correspondence determines a 
 birationality between the family of Pfaffian cubics and the family of 
 apolar cubics $F$. 
 On the other hand, computing the apolar quadrics to a Pfaffian cubic with 
 \cite{MAC} (cf. \ref{3.2} for apolarity), it can readily be checked that 
 there are in general no quadratic relations between them, 
 while the apolar quadrics to a cubic $F$ have nine quadratic relations: As 
 we shall see in the 
 next section, the apolar quadrics define the restriction to the $5$-space 
 $P=P_S=L_{S}^{\bot}$ of 
 the Cremona transformation defined by all quadrics through $G$. The  
 inverse Cremona 
 transformation is defined by quadrics again, and since $P$ has codimension 
 9 in $\Pn {14}$, 
 there are at least nine quadrics containing the image of $P$, i.e. at least 
 nine quadratic 
 relations between the apolar quadrics.
 \begin{problem}{\label{1.6}} Find an alternative description of the apolar 
 cubic 
 $4$-folds $F$.\end{problem}
 The main results of this paper are stated and proved in the third section. 
 In preparation and of 
 independent interest are some properties of the Cremona transformation 
 defined by the Pfaffians 
 of the matrix $M$, which we collect in the next section.
 \setcounter{subsection}{6} 
 \subsection{Notation}{\label{1.7}} The kernel of an element $\alpha\in 
 \wedge^2V^*$ is the subspace 
 $$\ker\alpha=\{v\in V |\; \alpha (u\wedge v)=0,\; \forall u\in V\}.$$ A 
 subspace $U\subset V$ 
 is Lagrangian with respect to $\alpha$ if $\alpha|_{U}\equiv 0$. 
 Similarly we define for $g\in \wedge^2V$ the kernel $\ker g\subset 
 V^{*}$ and Lagrangian subspaces with respect to $g$ in $V^{*}$. 
 For any $g\in \wedge^2V$, we denote 
 by $|g|$ the subspace ${(\ker g)}^{\bot}\subset V$, and call it the  
 support of $g$. For $g\in G$ the support $|g|$ is of course the rank 2 
 subspace 
 of $V$ represented by $g$. The rank of the support $|g|$ clearly 
 equals the rank of $g$. Similarly we define the support of 
 $\alpha$ in $V^{*}$.

 \section {Geometry of $\Gr (2,V)$ and its associated Cremona transformation}

 \subsection{\label{2.1}} In the Pl\"ucker coordinates $x_{ij}$, 
 the equations of $G=\Gr (2,V)$ are the $4\times 4$ Pfaffians of the  
 matrix $M$. 
 Let $q_{ij}=\Pfaff_{ij}M$ be the $4\times 4$ Pfaffians for $0\leq i<j\leq 5$ 
 and let $m=\Pfaff M$ be the $6\times 6$ Pfaffian of $M$. 
 While
 
 \begin{eqnarray*} 
 m&=&x_{05}x_{14}x_{23}-x_{04}x_{15}x_{23}-x_{05}x_{13}x_{24} \\ 
 &+&x_{03}x_{15}x_{24}+x_{04}x_{13}x_{25}-x_{03}x_{14}x_{25} \\ 
 &+&x_{05}x_{12}x_{34}-x_{02}x_{15}x_{34}+x_{01}x_{25}x_{34} \\ 
 &-&x_{04}x_{12}x_{35}+x_{02}x_{14}x_{35}-x_{01}x_{24}x_{35} \\ 
 &+&x_{03}x_{12}x_{45}-x_{02}x_{13}x_{45}+x_{01}x_{23}x_{45}, 
 \end{eqnarray*}
 the quadrics $q_{ij}$ are:
 \[\begin{array}{cc}
 q_{45}=x_{03}x_{12}-x_{02}x_{13}+x_{01}x_{23}& 
 q_{35}=x_{04}x_{12}-x_{02}x_{14}+x_{01}x_{24} \\ 
 q_{34}=x_{05}x_{12}-x_{02}x_{15}+x_{01}x_{25}& 
 q_{25}=x_{04}x_{13}-x_{03}x_{14}+x_{01}x_{34} \\ 
 q_{24}=x_{05}x_{13}-x_{03}x_{15}+x_{01}x_{35}& 
 q_{23}=x_{05}x_{14}-x_{04}x_{15}+x_{01}x_{45} \\ 
 q_{15}=x_{04}x_{23}-x_{03}x_{24}+x_{02}x_{34}& 
 q_{14}=x_{05}x_{23}-x_{03}x_{25}+x_{02}x_{35} \\ 
 q_{13}=x_{05}x_{24}-x_{04}x_{25}+x_{02}x_{45}& 
 q_{12}=x_{05}x_{34}-x_{04}x_{35}+x_{03}x_{45} \\ 
 q_{05}=x_{14}x_{23}-x_{13}x_{24}+x_{12}x_{34}& 
 q_{04}=x_{15}x_{23}-x_{13}x_{25}+x_{12}x_{35} \\ 
 q_{03}=x_{15}x_{24}-x_{14}x_{25}+x_{12}x_{45}& 
 q_{02}=x_{15}x_{34}-x_{14}x_{35}+x_{13}x_{45} \\ 
 q_{01}=x_{25}x_{34}-x_{24}x_{35}+x_{23}x_{45}& 
 \end{array} 
 \] 
 Notice that $(-1)^{i+j-1}q_{ij}$ is precisely the partial of $m$ with 
 respect to $x_{ij}$ i.e. 
 $$3m=\sum_{0\leq i<j\leq 5}(-1)^{i+j-1}x_{ij}q_{ij}.$$ 
 The Pfaffians $q_{ij}$, define a Cremona transformation (cf. 
 \cite{ES}) 
 $$\varphi :\Pn {14}\rto\Pn {14}.$$ In 
 fact $$q_{ij}(q_{st})=mx_{ij},$$ so the Cremona transformation is its own 
 inverse. The Cremona 
 transformation contracts precisely all secants to $G$. The exceptional 
 divisor lying over $G$ 
 in the Cremona transformation is mapped to a cubic hypersurface 
 $K^{\prime}$, the secant 
 variety of a variety $G^{\prime}$ which in turn is isomorphic to $G$.
 
 \subsection{\label{2.2}} 
 In terms of $2$-vectors $g\in \wedge^2V$ this Cremona transformation is the 
 composition 
 $$\varphi:\wedge^2V\to \wedge^4V\to \wedge^2V^{*},\quad g\mapsto g\wedge 
 g\mapsto 
 g\wedge g-,$$ where the last map is the natural isomorphism 
 (canonical up to scalars). 
 Therefore the target space of the Cremona transformation is naturally 
 identified with $\Pnd {14}$. 
 It is a morphism 
 on the complement of $G=\Gr (2,V)$, it is birational on the complement of 
 $K$, while 
 $K\setminus G$ is mapped to $G^{*}= \Gr (V,2)$. In fact, $g\wedge g=0$ for 
 $g\in 
 G$, while $\ker\varphi(g)=|g|$ when $g\in K\setminus G$.

 The fiber of $\varphi$ over a point 
 $\alpha\in G^{*}$ is 
 $$\varphi^{-1}(\alpha)=\{g\in 
 {\bf P}(\wedge^2 V)||g|=\ker\alpha\}={\bf P}(\wedge^2\ker\alpha),$$ a 
 $5$-dimensional space 
 which intersects $G$ in a quadric hypersurface $\Gr  (2,\ker\alpha)$.
 
 \subsection{\label{2.3}} The preimage under $\varphi$ 
 of a line in $G^*$ is a rational scroll ruled in $5$-dimensional spaces. 
 For this, 
 first note that the points where only the quadrics $q_{01}$ and $q_{02}$ 
 are nonzero, 
 are mapped to a line on $G^{*}$. On the other hand, by inspection, all the 
 quadrics 
 $q_{ij}$, except $q_{01}$ and $q_{02}$, vanish on the union of $G$ and the 
 cubic scroll defined 
 by the $2\times 2$ minors of 
 \[ 
 \left(\begin{array}{ccc} 
 x_{13}&x_{14}&x_{15}\\ x_{23}&x_{24}&x_{25}\end{array}\right) 
 \] 
 inside the $8$-dimensional space 
 $Z(x_{01},x_{02},x_{03},x_{04},x_{05},x_{12})$. 
 By homogeneity on $G^{*}$ the preimage 
 of any line is a 6-fold cubic scroll.
 
 \subsection{\label{2.4}} There are two kinds of planes in $G^*$. First we 
 have the planes representing all 
 lines in a plane, and second we have the planes representing all lines 
 through a point in a 
 $3$-space. The points where all the quadrics $q_{ij}$ except 
 $q_{01},q_{02}$ and $q_{12}$ 
 vanish, are all mapped to a plane in $G^*$ of the first kind. On the other 
 hand, by inspection, 
 all these quadrics vanish on the union of $G^*$ and a $7$-dimensional 
 subvariety of degree $6$ 
 defined by the $2\times 2$ minors of 
 \[\left(\begin{array}{ccc}x_{03}&x_{04}&x_{05}\\ x_{13}&x_{14}&x_{15}\\ 
 x_{23}&x_{24}&x_{25}\end{array}\right)
 \]
 inside the $11$-dimensional space $Z(x_{01},x_{02},x_{12})$. By 
 homogeneity the preimage 
 under $\varphi$ of any plane of the first kind is a $7$-fold of degree $6$ 
 in a $\Pn {11}$.
 The points where all the quadrics $q_{ij}$ except $q_{01},q_{02}$ and $q_{03}$ 
 vanish are all mapped to a plane in $G^*$ of the second kind. On the other 
 hand, by inspection, 
 all these quadrics vanish on the union of $G^*$ and the $7$-dimensional 
 complete intersection 
 of the two quadrics $q_{04}$ and $q_{05}$ inside the $9$-dimensional space 
 $Z(x_{01},x_{02},x_{03}, x_{04}, x_{05})$. By homogeneity the preimage 
 under $\varphi$ of any 
 plane of the second kind is a $7$-fold complete intersection 
 of two quadrics in a $\Pn {9}$.

 \subsection{\label{2.5}} Next we consider a tangent space to $G$. Without 
 loss of 
 generality we may consider the line spanned by 
 $L_{01}=<e_0,e_1>\subset {\bf P}(V)$ corresponding to the point 
 $p_{01}=(1,0,\ldots ,0)$ on 
 $G$. Let 
 \[ N_{01}=\left(\begin{array}{cccc}x_{02}&x_{03}&x_{04}&x_{05} \\ 
 x_{12}&x_{13}&x_{14}&x_{15}\end{array}\right) 
 \] 
 \setcounter{theorem}{5} 
 \begin{lemma}{\label{2.6}} $N_{01}$ has rank 1 on $G$ precisely at the 
 points which 
 correspond to lines which meet $L_{01}$. In fact the tangent 
 space to $G$ at $p_{01}$ is defined by $x_{ij}=0,\; 2\leq i<j\leq 5$, and 
 the $2\times 
 2$ minors of $N_{01}$ define the contact cone inside this tangent 
 space.\end{lemma} 
 \proof When 
 $x_{ij}=0$ for $2\leq i<j\leq 5$, then the Pl\"ucker quadrics reduce to the 
 minors of 
 $N_{01}$. On the other hand when this matrix has rank 1, i.e. $$\alpha 
 (x_{02},x_{03},x_{04},x_{05})+\beta 
 (x_{12},x_{13},x_{14},x_{15})=(0,0,0,0)$$ then it is the 
 Grassmannian point of the line $$(\beta e_0+\alpha e_1)\wedge 
 (x_{02}e_2+x_{03}e_3+x_{04}e_4+x_{05}e_5)$$ which is a general line which 
 meet $L_{01}$.\qed\vspace{.2 in}
 Next, we associate to each line in $K$ a rank $4$ 
 subspace of $V^{*}$ (cf. \cite{BD}).
 
 \begin{lemma}{\label{2.7}} 
 Let $l\in K$ be a line that does not intersect $G$, then there is a unique 
 $4$-space 
 $V^{*}(l)\subset V^{*}$, which is Lagrangian with respect to every element 
 of $l$.\end{lemma}
 \proof Let $g,g'$ be two points that span $l$. Then $\ker g$ and 
 $\ker g'$ are both $2$-dimensional subspaces of $V^{*}$ (rank$g$=rank$g'$=4). 
 Since $L$ does not intersect $G$, the two kernels span a $4$-space. We 
 may choose a 
 basis $<e^{*}_{0},\ldots, e^{*}_{5}>$ for $V^{*}$ such that $<e^{*}_{0}, 
 e^{*}_{1}>$ and $<e^{*}_{4}, e^{*}_{5}>$ are the 
 two kernels. Thus the skewsymmetric matrices $M(g)$ and $M(g')$ have 
 coefficients 
 $$M(g)=(a_{ij}),\; a_{ii}=0, a_{ij}=-a_{ji}, a_{0j}=a_{1j}=0$$ 
 and 
 $$M(g')=(b_{ij}),\; b_{ii}=0, b_{ij}=-b_{ji}, b_{4j}=b_{5j}=0.$$ 
 Let $h=sg+tg'$ be a general point on $l$, then 
 $M(h)=(sa_{ij}+tb_{ij})$. Now $l\subset K$, so $M(h)$ has rank 
 $4$ (independant of $s$ and $t$). Therefore 
 $\Pfaff M(h)=3st(ta_{45}\Pfaff_{45}M(g') + sb_{01}\Pfaff_{01}M(g))=0$ for 
 every $s$ and $t$. But 
 since $M(g)$ and $M(g')$ have rank $4$, both $\Pfaff_{01}M(g)$ and 
 $\Pfaff_{45}M(g')$ are 
 nonzero, so $a_{45}=b_{01}=0$. This means that that the $4$-space 
 $<e_{0},e_{1},e_{4},e_{5}>$ is 
 Lagrangian for all $h\in l$. On the other hand any common  Lagrangian
 $4$-space has to contain $\ker g$ and $\ker g'$ so it is 
 unique.\qed\vspace{.2 in}

 Finally we investigate certain subvarieties associated to secant lines to 
 $G^*$. 
 For a proper secant line $l^{*}$, we associate a quadric 
 surface $Q(l^{*})\subset G^{*}$. If the secant line $l^{*}$ intersects 
 $G^{*}$ in two points, the quadric surface $Q(l^{*})$ parametrizes all lines 
 that intersect the lines represented by these two points. Notice that 
 the support $|\alpha|$ 
 is common for every rank $4$ element $\alpha\in l^{*}$, we denote this by 
 $U(l^{*})$. Thus 
 $Q(l^{*})\subset {\Gr }(2,U(l^{*}))\subset {{\bf P}}(\wedge^2U(l^{*}))$. 
 Since 
 $U(l^{*})$ is $4$-dimensional, $\Gr 
 (2,U(l^{*}))$ is a smooth quadric hypersurface. Polarity with respect to 
 this quadric defines a 
 correlation $p_U:{\bf P}(\wedge^2U(l^{*}))\to {\bf P}(\wedge^2U(l^{*})^*)$. 
 Now $l^*\subset {\bf 
 P}(\wedge^2U(l^{*}))$, so $p_U(l^*)^{\bot}\subset {\bf 
 P}(\wedge^2U(l^{*}))$ is a $3$-space $P(l^*)$, the span of $Q(l^*)$.\par 
 When $l^*\subset {\bf P}(\wedge^2U)$ for some $4$-dimensional subspace 
 $U\subset V^*$, we 
 denote by $P_U(l^*)$ the $3$-space $p_U(l^*)^{\bot}\subset {\bf 
 P}(\wedge^2U)$ polar to $l^{*}$ with respect to the quadric $\Gr (2,U)$. 
 The pairs $(l^{*},Q(l^{*}))$ therefore fit into the following 
 incidence:
 $$I_Q=\{(l^{*},Q)|\exists U\subset V^*, \dim 
 U=4, l^{*}\subset {\bf P}(\wedge^2U), Q=P_U(l^*)\cap G^{*}\}$$
 
 \begin{lemma}{\label{2.8}} Let $(l^{*},Q)\in I_Q$. Then $l^{*}$ 
 is a secant line to $G^{*}$, a tangent line to $G^{*}$ or is contained in 
 $G^{*}$. 
 It is a secant line if and only if $Q$ is smooth. It is a tangent 
 line if and only if $Q$ has rank $3$. It is contained in $G^{*}$ if and 
 only if $Q$ has rank $2$. The projection onto the second factor is 
 onto the variety of quadric surfaces in $G^{*}$, whose span is not 
 contained in $G^{*}$. The first projection is one to one except over the 
 lines in $G^{*}$, where each fiber is a $\Pn 2$.\end{lemma}
 \proof The definition above of $Q(l^{*})$ clearly extends to tangent 
 lines: If $ l^{*} $ is a tangent line to $G^{*}$, then $|\alpha|=U$  
 is constant for every rank $4$ element $\alpha\in l^{*}$, and 
 $Q(l^*)$ is a singular $\Pn 3$ section of $\Gr (2,U)$. If $ l^{*} $  
 is a line in $G^{*}$, then $\cup_{\alpha}|\alpha|$ for $\alpha\in 
 l^{*}$ is a $3$-dimensional subspace $U_{0}\subset V^{*}$. Any 
 $4$-space $U\supset U_{0}$, gives rise to a quadric $Q\subset G^{*}$ of 
 rank $2$. One of the planes of $Q$ is common for all $U\supset 
 U_{0}$, it is the plane $\Gr (2,U_{0})$. In all cases it is 
 straightforward to check that the 
 span of the quadric surface is polar to $l^*$ with respect to the quadric 
 $\Gr (2,U)$.
 Thus we are left to prove 
 that the second projection is onto the family of quadric surfaces, 
 whose span is not contained in $G^{*}$.
 
 Let $Q$ be a quadric surface of this kind. First note that the span $P_Q$ 
 of $Q$ is contained in 
 $K^*$. The restriction of the inverse Cremona transformation to $P_Q$, is 
 clearly the 
 constant map onto some point on $G$. Thus $P_Q$ is contained in fiber of 
 this inverse 
 transformation. By \ref{2.2}, the whole fiber is a ${\bf P}(\wedge^2U)$ 
 which intersect $G^*$ in the 
 rank $6$ quadric $\Gr (2,U)$. Therefore the intersection of this quadric 
 with a $3$-space has 
 rank at least 2. Furthermore $P_Q=P_U(l^*)\subset {\bf P}(\wedge^2U)$ for 
 some 
 $l^*\subset {\bf P}(\wedge^2U)$. The lemma now follows.\qed\vspace{.2 in}
 
 \begin{lemma}{\label{2.9}} Let $Q$ be a quadric 
 surface in $G^*$, whose span is not contained in $G^*$, and let $l^*$ be a 
 line in $G^*$ such 
 that $(l^*,Q)$ belong to the incidence $I_Q$. Then the preimage 
 $\varphi^{-1}(Q)$ under 
 the Cremona transformation is a generic variety, the intersection of two 
 rational quartic 
 scrolls in a rank 4 quadric hypersurface inside $(l^*)^{\bot}$. 
 In particular, it has degree $10$ and 
 codimension $5$ inside $(l^*)^{\bot}$.\end{lemma} 
 \proof 
 For this, consider first the secant line $l^*$ of $G^*$ spanned by the 
 points with coordinates 
 $y_{st}=0$ for $(st)\neq (01)$ and $(st)\neq (23)$ respectively. Then 
 $Q(l^*)=Z(y_{02}y_{13}-y_{12}y_{03})$ inside the 
 $3$-space $P(l^*)=Z(y_{st}|(st)\notin\{02,03,12,13\})$.
 Therefore the preimage of $Q(l^*)$ under the Cremona transformation is 
 defined by the 
 $11$ quadrics 
 $$V(01,23)=<q_{st}|(st)\notin\{02,03,12,13\}>.$$ 
 These quadrics vanish on the union of $G^*$ and 
 $$Z(01,23)=Z(V(01,23))\cap Z(x_{01}, x_{23}).$$ 
 By inspection, 
 $Z(01,23)$ 
 is the locus inside the hyperplanes $Z(x_{01})$ and $Z(x_{23})$ 
 where the matrices 
 \[ N_{01}=\left(\begin{array}{cccc}x_{02}&x_{03}&x_{04}&x_{05}\\ 
 x_{12}&x_{13}&x_{14}&x_{15}\end{array}\right) 
 \; {\rm and} \; N_{23}=\left(\begin{array}{cccc}x_{02}&x_{12}&x_{24}&x_{25}\\ 
 x_{03}&x_{13}&x_{34}&x_{35}\end{array}\right) 
 \] 
 both drop rank. 
 The invariants of $Z(01,23)$ are easily computed in 
 MACAULAY \cite{MAC}. The matrices $N_{01}$ 
 and $N_{23}$ each drop rank on a rational quartic scroll of codimension 3. 
 The two 
 matrices have the $2\times 2$ submatrix 
 $\left(\begin{array}{cc}x_{02}&x_{03}\\ x_{12}&x_{13}\end{array}\right)$ in 
 common, so $Z(01,23)$ is the intersection of two rational normal quartic 
 scrolls inside a rank 4 
 quadric. For the degree we note that this subvariety is a degeneration of 
 the intersection 
 of two codimension 2 cycles of bidegree $(1,3)$ on a rank $6$ quadric. Thus 
 the degree is 
 10.\par
 By homogeneity \ref{2.9} follows for smooth quadric surfaces. The 
 computations for quadric 
 surfaces of rank 3 and 2 are quite similar and straightforward to 
 check.\qed\vspace{.2 in}

 \section{The variety of sums of powers ${{\rm VSP}(F,10)}$}
 \subsection{\label{3.1}} For a homogeneous polynomial $f$ of degree $d$ in $n+1$ 
 variables, 
 which define the hypersurface $F=Z(f)\subset \Pn n$, we define the variety 
 of sums of powers 
 as the closure 
 $$VSP(F,s) = \overline{\{\{<l_1>,\ldots,<l_s>\}\in Hilb_s(\check {\Pn n})\mid\exists\lambda_i \in \C: f=\lambda_1l_1^d+\ldots+\lambda_s l_s^d\}}$$ 
 of the set of powersums presenting $F'$ in the Hilbert scheme (cf. 
 \cite{RS}). We study these powersums 
 using apolarity.
 
 \subsection{Apolarity}{\label{3.2}}(cf. \cite{RS}). 
 Consider $R=\C[x_0,\ldots,x_n]$ and 
 $T=\C[\partial_0,\ldots,\partial_n]$. $T$ acts on $R$ by 
 differentiation: 
 $$\partial^{\alpha}(x^{\beta}) = \alpha!{\binom{\beta}{\alpha} } 
 x^{\beta-\alpha}$$ 
 if $\beta \geq \alpha$ and 0 otherwise. Here $\alpha$ and $\beta$ are 
 multi-indices, 
 ${\binom{\beta}{\alpha}} = \prod {\binom{\beta_i}{\alpha_i}}$ and so on. 
 In particular we have a perfect pairing between forms of degree d and 
 homogeneous differential operators of order d. 
 Note that the polar of a form $f \in R$ in a point $a \in \Pn n$ is  
 given by $P_a(f)$ for 
 $a=(a_0,\ldots,a_n) $ and $P_a = \sum a_i \partial_i \in T$.
 One can interchange the role of R and T by defining 
 $$x^{\beta}(\partial^{\alpha}) = \beta!{\binom{\beta}{\alpha} } 
 \partial^{\alpha-\beta}.$$ 
 With this notation we have for forms of degree n 
 $$P_a^d(f)=f(P_a^d)=d!f(a).$$ 
 Moreover $$f(P_a^m)=0 \iff f(a)=0 \eqno(*)$$ if $m \geq d$. 
 More generally we say that homogeneous forms $f\in R$ and $D\in T$ 
 are {\bf apolar} if $f(D)=D(f)=0$ (According to 
 Salmon (1885) \cite{Sal} the term was coined by Reye).
 
 Apolarity allows us to define Artinian Gorenstein 
 graded quotient rings of $T$ 
 via forms: For $f$ a homogeneous form of degree $d$ and $F =Z(f) \subset 
 \Pn n$ define 
 $$ F^{\bot} = f^{\bot} = \{D \in T | D(f)=0 \}$$ 
 and $$A^F = T/F^{\bot}.$$
 An obvious but useful identity in apolarity is 
 $$(f^{\bot}:D)=D(f)^{\bot}  $$ 
 for any homogeneous $D \in T$.
 The socle of $A^F$ is in degree d. Indeed 
 $P_a(D(f)) = 0 \hskip3pt \forall P_a \in T_1 \iff D(f) = 0 \hskip3pt 
 or \hskip3pt D \in T_d$. 
 In particular the socle of $A^F$ is 1-dimensional, and $A^F$ is 
 indeed Gorenstein and is called the {\bf apolar Artinian 
 Gorenstein ring} of $F\subset \Pn n$.
 Conversely for a graded Gorenstein ring $A = T/I$ with socledegree $d$, 
 multiplication in A induces a linear form $f\colon {\rm Sym}_d(T_1) \to \C$ 
 which can be identified with a homogeneous 
 polynomial $f \in R$ of degree $d$. This proves: 
 \setcounter{theorem}{2} 
 \begin{lemma}{\label{3.3}} {\rm (Macaulay, \cite{Mac})} 
 The map $F \mapsto A^F$ is a bijection between hypersurfaces $F=V(f) 
 \subset \Pn n$ of 
 degree d and graded Artinian Gorenstein quotient rings $A = T/I$ of T 
 with socledegree d.\end{lemma}
 In the following we identify R with homogeneous 
 coordinate ring of $\Pn n$ 
 and T with the homogeneous coordinate ring of the dual space $ 
 \Pnd n$. Let $F=Z(f) \subset \Pn n$ 
 denote a hypersurface of degree d. We call a subscheme 
 $\Gamma \subset \check {\Pn n}$ 
 $\bf {apolar} $ to $F$, if the homogeneous ideal $I_{\Gamma} \subset 
 F^{\bot} \subset T$.
 \begin{lemma}{\label{3.4}} Let $l_1,\ldots, l_s$ be linear forms in $R$, 
 and let $L_i \in \Pnd n$ be 
 the corresponding points in the dual space. Then 
 $f = \lambda_{1}l_1^d+\ldots+\lambda_{s}l_s^d$ for some 
 $\lambda_{i}\in \C^{*}$ if and only if $\Gamma = \{ L_1,\ldots 
 ,L_s\} \subset \check {\Pn n}$ is apolar to $F=Z(f)$.\end{lemma}
 \proof 
 Assume $f = \lambda_{1}l_1^d+\ldots+\lambda_{s}l_s^d$. If $g\in 
 I_{\Gamma}$, then $g(l_i^d)$=0 for all $i$ by (*), so by 
 linearity $g\in F^{\bot}$. Therefore $\Gamma$ is apolar to $F$. \par 
 For the converse, assume that $I_{\Gamma}\subset F^{\bot}$. Then we have 
 surjective maps between 
 the corresponding homogeneous coordinate rings 
 $$T\to A_{\Gamma}=T/I_{\Gamma}\to A^F.$$ 
 Consider the dual inclusions of the degree $d$ part of these rings:  
 $$\Hom (A^F_d,\C)\to \Hom ((A_{\Gamma})_d,\C)\to \Hom (T_d,\C).$$ 
 $D\mapsto D(f)$ generates the first of these spaces, while the second is 
 spanned by the forms 
 $D\mapsto D(l_i^d)$. Thus $F'$ lies in the span of the $l_i^d$. 
 \qed\vspace{.2 in} 
 \setcounter{subsection}{4} 
 \subsection{\label{3.5}} We return now to the case of cubic $4$-folds. 
 In the notation of \ref{1.3} we are given 
 a general $5$-dimensional space $P_S\subset\Pn {14}$ defining the apolar 
 Artinian Gorenstein ring 
 $A^{F}$ of some cubic $F=Z(f)\subset \check P_S$. By \cite{AH} (cf.  
 \cite[1.1]{RS}) the minimal degree 
 of finite apolar subschemes of a general cubic $4$-fold is $10$. We first 
 show that $F$ is 
 general in this sense. 
 Following \ref{3.4} we study ideals of finite subschemes contained 
 in $(F)^{\bot}$. Now, 
 $(f)^{\bot}$ is generated by 15 quadrics. On the one hand these are  
 precisely the orthogonal 
 complement in $T_2$ of the partials of $F$ in $R_2$. On the other hand 
 they are precisely the 
 quadrics defining $G$ restricted to $P_S$. Consider the morphism 
 $\varphi_P$ defined by these 
 quadrics generating $(f)^{\bot}$. On the one hand it is the composition of 
 the $2$-uple 
 embedding and the projection from the partials of $f$, on the other hand it 
 is the restriction 
 of the Cremona transformation $\varphi$ to $P_S$.

 \setcounter{theorem}{5} 
 \begin{lemma}{\label{3.6}} The minimal degree of a finite apolar subscheme 
 of $F$ is 10.\end{lemma} 
 \proof 
 A subscheme of length $9$ lies on at least $12$ quadrics. Therefore  
 an apolar subscheme of length $9$ would be mapped by $\varphi$ to a scheme 
 of length $9$ in a 
 plane. If the inverse Cremona transformation restricted to this plane is 
 birational, then 
 the image of this plane in $\Pn {14}$ would have degree at most $4$, so its 
 intersection with $P_S$ could not be of length 9. 
 Therefore the plane is contained in $K^*$. But then the subscheme of 
 length $9$ must be contained in $G^*$. If the plane intersects $G^*$ in a 
 conic section, 
 then the subscheme of length $9$ is contained in this conic section, i.e. 
 in a line in $P$. This 
 is absurd, since any apolar subscheme must span $P$. Thus the plane must 
 be contained in $G^*$. 
 There are two kinds of planes in $G^*$. 
 We know from \ref{2.4} 
 that the preimage is of codimension $4$ and $2$, respectively, in their linear span.
 Again, since apolar subschemes span $P$, this is impossible. \qed\vspace{.2 in}

 A length $10$ subscheme $\Gamma$ lies on 
 at least $11$ quadrics. Therefore if $\Gamma$ is apolar to $F$, the image 
 $\varphi_P(\Gamma)$ 
 spans at most a $\Pn 3$. 
 Hence it is natural to identify a powersum of length $10$ presenting $f$ 
 with some $10$-secant $\Pn 3$ to $\varphi(P_S)$. Let 
 $VSP_G(F,10)$ be the closure in the Grassmannian 
 ${\bf G} (4,\wedge^2V^{*})$ of $10$-secant $3$-spaces 
 to $\varphi(P_S)$. 
 We call this a 
 Grassmannian compactification of the set of powersums of $f$. 
 We proceed in a few lemmas to study $VSP_G(F,10)$, and in fact to show 
 that every $\Pn 3$ of this 
 family intersect $\varphi(P_S)$ in a finite subscheme of length $10$. 
 Therefore the Grassmannian 
 compactification coincides with the Hilbert scheme compactification. 
 We now make all this more precise.

 Let $VSP_G(F,10)$ be the closure in the Grassmannian 
 ${\bf G} (4,\wedge^2V^{*})$ of $3$-spaces that intersect 
 $\varphi(P_S)$ in a finite subscheme of length $10$ apolar to 
 $F$.

 \begin{theorem}{\label{3.7}} $VSP_G(F,10)$ is isomorphic to the family of 
 secant lines to $G^*\cap 
 L_{S}$, i.e. to ${\rm Hilb}_2(S)$ where $S$ is the $K3$ surface $S=G^*\cap 
 L_{S}$.\end{theorem} 
 \proof We first define an injective map $\rho: {\rm Hilb}_2(S)\to 
 VSP_G(F,10)$. 
 For this let $l^{*}$ be a secant or tangent line to $S$. It is defined by 
 its intersection of length $2$ with $S$. Since $S$ is general, it 
 contains no lines, and is the intersection of quadrics, so the 
 corresponding point in ${\rm Hilb}_2(S)$ is unique. Now, let 
 $Q(l^{*})\subset G^{*}$ be the quadric surface associated to 
 $l^{*}$. According to \ref{2.9} the preimage $\varphi^{-1}( Q(l^{*}))$ is 
 $7$-dimensional of degree $10$ and 
 contained in $(l^{*})^\bot$. But $l^{*}\subset P_S^\bot=L_S$, means that 
 $P_S\subset (l^{*})^\bot$, so $\Gamma_{l^{*}}=P_S\cap \varphi^{-1} 
 (Q(l^{*}))$ 
 is an apolar 
 subscheme to $F$ of length $10$ as long as the intersection is proper. 
 Clearly this is so for the general $P_S$ and secant line $l^{*}$. 
 Furthermore, the image $\varphi( \Gamma_{l^{*}})$ is always the 
 intersection of the $3$-space $P(l^{*})=<Q(l^{*})>$ with $\varphi(P_S)$, 
 so the 
 map $\rho$ is welldefined. Injectivity follows from the $1:1$ 
 correspondance $l^{*}\leftrightarrow P(l^*)$.
 To show that there is a welldefined converse map we need a few lemmas. 
 Let $\Gamma\subset P_S$ be an apolar subscheme to $F$, such that 
 $\varphi(\Gamma)$ spans a $\Pn 3$, which we denote by $P_{\Gamma}$, and 
 assume that $P_{\Gamma}\in 
 VSP_G(F,10)$. Thus $P_{\Gamma}$ lies on the closure of the variety 
 of $10$-secant $3$-spaces to $\varphi(P_S)$. 
 In this notation: 
 \begin{lemma}{\label{3.8}} Let $P_\Gamma \in VSP_G(F,10)$, then 
 $P_{\Gamma}=P(l^{*})$ for some secant 
 or tangent line $l^{*}$ to $K^{*}$.\end{lemma} 
 \proof This proof depends on two lemmas. The first one is 
 interesting on its own: 
 \begin{lemma}{\label{3.9}} Let $P_{\Gamma} \in VSP_G(F,10)$ as above, then 
 $\Gamma\subset K\cap P_S$, 
 and $P_{\Gamma}\subset K^*$.\end{lemma} 
 \proof 
 Let $P_{\Gamma}\in 
 VSP_G(F,10)$. We may assume first that $\Gamma$ is smooth of length  
 $10$ and spans $P_S$. Then $\Gamma\cap 
 G\subset P_S\cap G=\emptyset$, so by the Cremona transformation any point 
 of $\Gamma$ is mapped 
 to $G^*$ or to the complement of $K^*$. The image $\varphi (\Gamma)$ spans 
 a $3$-space 
 $P_{\Gamma}$. The inverse Cremona restricted to $P_{\Gamma}$ is defined by 
 the quadrics 
 through $G^*\cap P_{\Gamma}$. Assume that $K^*\cap P_{\Gamma}\not= 
 P_{\Gamma}$. Then the 
 restriction of the inverse Cremona transformation to $P_{\Gamma}$ is 
 birational onto its 
 image. Since the image intersects $P_S$ in $\Gamma$ which in turn span 
 $P_S$, the image must span 
 at least an $8$-space. But in that case the degree of the image is 7 or 8 
 so it cannot 
 contain $\Gamma$. Therefore $P_{\Gamma}\subset K^*$, and $\Gamma\subset K$. 
 Now, 
 $VSP_G(F,10)$ is a compactification of a set of $10$-secant $\Pn 3$'s to 
 $\varphi(P_S)$. 
 Therefore any $P_{\Gamma}$ in the closure must also be contained in $K^*$. 
 Furthermore the intersection $\Gamma=P_{\Gamma}\cap \varphi(P_S)$ must 
 be of length at least $10$. \qed\vspace{.2 in}

 For \ref{3.8} we note that $\Gamma $ is mapped to $G^*$ by the Cremona 
 transformation. In 
 the above notation
 \begin{lemma}{\label{3.10}} If $P_{\Gamma}\in VSP_G(F,10)$, then 
 $P_\Gamma\cap 
 G^{*}$ is a quadric 
 surface.\end{lemma} 
 \proof Consider the inverse Cremona transformation. Since 
 $P_{\Gamma}\subset K^*$, the image of $P_{\Gamma}\setminus G^{*}$ is 
 contained in $G$. 
 By \ref{2.2}, any fiber over $G$ of the inverse Cremona is 
 a $\Pn 5$ which intersects $G^{*}$ in a quadric $4$-fold. Therefore the 
 $3$-space $P_\Gamma$ 
 meets each fiber of this Cremona transformation in a linear space, while 
 each of these 
 linear spaces again intersect $G^{*}$ in at least a quadric. In the proof 
 of \ref{3.9} we saw that the 
 restriction of the inverse Cremona transformation to $P_\Gamma$ is not 
 birational. Therefore each 
 fiber is a line, plane or all of $P_\Gamma$.
 Assume first that $P_\Gamma$ is contained in $G^*$. Then the inverse image 
 by the Cremona transformation 
 is a quadric hypersurface in a $\Pn 9$ which intersect $G$ in a 
 Grassmannian $G(2,5)$. But 
 $P_S$ does not intersect $G$, so $P_S$ intersects the $\Pn 9$ in at most a 
 plane. In particular, 
 $\Gamma$ is contained in a plane. But this would mean that $\Gamma$  
 is an apolar subscheme of $F$ contained in a plane, clearly absurd.
 Therefore $P_\Gamma$ is not all 
 contained in $G^*$. Thus $P_\Gamma\cap G^*$ is a subvariety defined by 
 quadrics, and each fiber of the inverse Cremona transformation restricted to 
 $P_{\Gamma}$ is a line or a plane.
 When the general fiber is a plane, this plane must intersect $G^*$ in at 
 least a conic section, so the intersection $G^*\cap P_{\Gamma}$ must 
 therefore be a quadric 
 surface or the union of a plane and a line. On the other hand it is not 
 hard to check that 
 the latter is never the intersection of a $3$-space with $G^*$, so if the 
 fibers are planes 
 then the intersection is a quadric surface.
 When the fibers are lines, then 
 $G^*\cap P_{\Gamma}$ is a quadric surface, the union of a plane and a line 
 or a curve 
 defined by quadrics with one secant line through each general point of 
 $P_\Gamma$, i.e. a twisted cubic curve or the union of two lines. 
 The union of a 
 plane and a line never 
 occurs 
 as the intersection of a $3$-space with $G^*$, while the other cases 
 do occur as 
 intersections with the Grassmannian. Note, that there could be no extra 
 points of 
 intersection in addition to these curves, since these would lie on proper 
 trisecants to the 
 Grassmannian, which is absurd since the Grassmannian is cut out by 
 quadrics. \par 
 Now, the preimage under $\varphi $ of a line is a 6-fold cubic 
 scroll by \ref{2.3}, so the preimage of a conic section must be a 6-fold 
 scroll of degree 6 and the 
 preimage of a twisted cubic curve must be a 6-fold scroll of degree 9. But 
 $\Gamma $ is the 
 intersection of the $5$-space $P_S$ with the preimage under the Cremona 
 transformation of 
 $G^*\cap P_{\Gamma}$. If $\Gamma$ is 0-dimensional, the length of this 
 intersection cannot 
 exceed the degree of the corresponding preimages. Since these all have 
 degree less than 10, 
 these cases are excluded and we conclude that $G^*\cap P_{\Gamma}$ is 
 a quadric surface or all of $P_{\Gamma}$. If 
 $\Gamma$ has positive dimension, then, since $\Gamma$ lies in the closure 
 of finite schemes, 
 $G^{*}\cap P_{\Gamma}$ is a quadric surface or all of $P_\Gamma$. 
 Since the latter is already excluded, the lemma follows. \qed\vspace{.2 in}
 
 We denote the quadric surface $P_{\Gamma}\cap G^{*}$ by $Q_{\Gamma}$. 
 To conclude the proof of \ref{3.8} it suffices to show that $Q_{\Gamma}$ is 
 irreducible, but by \ref{2.8} the quadric $Q_\Gamma$ is reducible only if 
 $l^*$ is contained in $G^*$. 
 This is excluded by our generality assumption of $S$. \qed\vspace{.2 
 in}
 
 For \ref{3.7} it remains to note that $\Gamma$ is apolar to $F$ only if 
 $P_S\subset(l^{*})^{\bot}$, where $P(l^{*})=P_{\Gamma}$. \qed\vspace{.2 in}
 
 \begin{proposition}{\label{3.11}} For a general $K3$ surface $S$, let $F(S)$ 
 be the associated apolar 
 cubic $4$-fold. Then the $4$-folds $VSP(F(S),10)$ and $VSP_G(F(S),10)$ 
 are isomorphic.\end{proposition} 
 \proof 
 To prove this we have to show that 
 for the 
 general $K3$ surface $S = G^* \cap L_S$, all the apolar subschemes 
 $\Gamma$ such that the span $P_\Gamma$ of $\varphi(\Gamma)$ is a $3$-space 
 contained in $VSP_G(F,10)$, 
 are finite. For this we shall see first the following: 
 \begin{lemma}{\label{3.12}} Let $P_\Gamma\in VSP_G(F,10)$, and assume that 
 $C$ is an integral 
 curve contained in $\Gamma$, then $C$ is a line.\end{lemma}
 $Proof.$ First, note that the restriction $\varphi$ to $P_S$ is an 
 embedding defined by a linear 
 system in $|2h|$, where $h$ embeds $C$ in $P_S$. Now, $\varphi(C)$ is 
 contained in a quadric 
 surface by \ref{3.10}. If $\varphi(C)$ lies in a plane, then clearly $C$ 
 must be a line. If 
 $\varphi(C)$ lies on a quadric cone (or a smooth quadric surface) $Q_{C}$, 
 then 
 it meets every ruling (or every line in one ruling respectively) in at 
 least $d$ points, 
 where $d$ is the degree of $C$ in $P_S$.
 Consider the preimage of a ruling of a quadric on $G$ under the 
 Cremona transformation. In \ref{2.3} we showed that the preimage of a line 
 $l$ is a cubic 6-fold scroll in a $\Pn 8_{l}$. When the quadric 
 $Q_{C}$ is 
 smooth, the intersection 
 $\cap_{l}\Pn 8_{l}$ over all lines $l$ in the ruling $r$ is a $4$-space 
 which we denote by $\Pn 4_{r}$. 
 Furthermore $\Pn 
 4_{r}\cap G$ is a rank $4$ quadric $3$-fold. When the quadric $Q_{C}$ is a 
 cone, the intersection $\cap_{l}\Pn 8_{l}$ is the $5$-space $\Pn 
 5_{r}=\varphi^{-1}(v)$, where $v={\rm Sing} Q_{C}$ is the vertex of  
 $Q_{C}$, and the intersection $\Pn 
 5_{r}\cap G$ is a smooth quadric. Now, if $\varphi(C)$ intersects each 
 line $l$ in the ruling $r$ in $a$ points, then $C\cap\Pn 8_{l}$ 
 contains at least $a$ points. If $a>d$, this means that $C\subset \Pn 
 8_{l}$ for each line $l$ in the ruling. So $C$ is contained in the 
 intersection $\Pn 4_{r}$ or $\Pn 5_{r}$. In either case $C$ must 
 intersect $G$, a contradiction.
 Therefore $\varphi(C)$ meets every ruling in at most $d$ points, i.e. 
 it must be a complete intersection of the quadric surface with 
 a surface of degree $d$ and $g(C)=(d-1)^2$. If $C$ is a plane curve, it is clearly a line.Ê Otherwise, Castelnuovo's bound for space curves applies and 
 $g(C)=(d-1)^2\leq \frac{1}{4}d^2-d+1=(\frac{d}{2}-1)^2$, which is possible 
 only if $d=1$. 
 \qed\vspace{.2 in}

 Let ${\cF}(K)$ be the 22-fold of lines on 
 the Pfaffian cubic 13-fold $K \subset {\bf P}^{14}$, 
 let $G(6,15)$ be the Grassmannian of 5-dimensional 
 projective subspaces $P_S \subset {\bf P}^{14}$. 
 Any line $l\in{\cF}(K)$, is mapped to a conic section $\varphi (l)\subset {\Pnd 
 {14}}$ by the Cremona transformation.
 On the other hand, by \ref{2.7}, there is a unique common 
 Lagrangian $4$-space $V^{*}(l)\subset V^*$ for all $g\in l$. When $l\subset 
 P_S\cap K$, then ${\bf P}(\wedge^2V^{*}(l))$ intersects $L_{S}=P_S^{\bot}$ 
 in a line $l^*$, which is neccessarily a secant line to $G^{*}$, therefore 
 also contained in $K^{*}$. In fact, let $P(l^*)\subset {\bf 
 P}(\wedge^2V^{*}(l))$ be the $3$-space associated to the secant line 
 $l^{*}$ to $K^{*}$ like in \ref{2.8} (i.e. the polar $\Pn 3$ to $l^*$ with 
 respect to the quadric $\Gr (2,V^{*}(l))$). It follows from \ref{2.9} that the 
 preimage of $P(l^*)$ under the Cremona transformation is contained in the 
 complement ${l^{*}}^{\bot}\subset \Pn {14}$. 
 But $l^{*}\subset L_S=P_S^{\bot}$, so $P_S\subset 
 {l^{*}}^{\bot}$. If the intersection between $P_S$ and the fiber is  
 proper inside ${l^{*}}^{\bot}$, then it is a finite subscheme of 
 length $10$. In other words $P(l^*)$ 
 is a 
 $10$-secant $\Pn 3$ to $\varphi (P_S)$ for general $l\subset P_S$.
 
 It remains to show that for general $P_S$ the intersection is always 
 proper. 
 Now clearly $\varphi (l)\subset {\bf P}(\wedge^2V^{*}(l))$, since 
 $g\in l$ implies that $V^{*}(l)$ is a Langrangian $4$-space with 
 respect to $g$, i.e. $\ker g\subset V^{*}(l)$. Therefore $l\subset 
 \Gamma$ only if $\varphi (l)\subset P(l^*)\subset{\bf 
 P}(\wedge^2V^{*}(l)).$

 Consider the incidence 
 \[ 
 I \ = \ 
 \{ (l,P_S)| 
 l\subset P_S,\; \varphi(l)\subset P(l^*) 
 \}, 
 \] 
 and let 
 $p: I \rightarrow {\cF}(K)$ 
 and 
 $q: I \rightarrow \Gr (6,\wedge^2V)$ 
 be the natural projections. The incidence naturally parametrizes the bad 
 locus, i.e. the set of 
 spaces $P_S$ such that there exist some $P_\Gamma\in VSP_G(F,10)$ where 
 $\Gamma$ contains a line. 
 Therefore it is enough to show that the second projection is not surjective. 
 In fact we will show that the codimension of $q(I)$ in 
 $\Gr(6,\wedge^2V)$ is at least $2$.
 
 Let $l \subset K$ be a general line; 
 in particular $l \subset K\setminus G$, 
 Then 
 \[ 
 p^{-1}(l) 
 = 
 \{ P_S = {\bf P}^5 \subset {\bf P}^{14}| 
 l \subset P_S, 
 \; 
 \varphi(l)\subset P(l^*) \}.\] 
 But $\varphi(l)$ and $P(l^*)$ are both subvarieties of ${\bf 
 P}(\wedge^2V^{*}(l))$. The cycle of $4$-spaces that contain the 
 conic section $\varphi(l)$ is a plane inside $G(4,V^{*}(l))$, so it is of 
 codimension $6$ in the 
 Grassmannian and the fiber 
 $p^{-1}(l)$ is of codimension $6$ among all $P_S$ that contain $l$, i.e. 
 the fiber have 
 dimension $30$. On the other hand the family of lines in the Pfaffian 
 cubic $K$ is of 
 dimension 22, while $\Gr (6,\wedge^2V)$ has dimension 54. Therefore the image of 
 $q$ has codimension 
 at least $2$.\qed\vspace{.2 in}

 \begin{corollary}{\label{3.13}} Let $P_S$ be general, then $VSP(F(S),10)$ is 
 isomorphic to the family of 
 secant lines to $G^*\cap L_{S}$, i.e. to ${\rm Hilb}_2(S)$ where $S$ is the 
 $K3$ surface 
 $S=G^*\cap L_{S}$. Furthermore, $VSP(F(S),10)$ is isomorphic to the  
 Fano variety of lines ${\cF}(F')$ of the Pfaffian 
 cubic fourfold $F'=F'(S)=P_S\cap K$, where $P_S=L_S^{\bot}$. 
 \end{corollary} 
 \proof The first statement follows directly from \ref{3.7} and 
 \ref{3.11}.
 The argument following the proof of \ref{3.12} shows that ${\rm 
 Hilb}_2(S)$ is isomorphic to the Fano variety of lines in $F'(S)$, as was 
 noticed by Beauville and 
 Donagi (cf. \cite{BD}).\qed\vspace{.2 in}
 \setcounter{subsection}{13} 
 \subsection{\label{3.14}} The corollary suggests an incidence correspondence 
 $$I=\{(g,\Gamma) | g\in \Gamma\}\subset F'(S)\times 
 VSP(F(S),10)\cong F'\times{\cF}(F').$$

 The second projection is clearly $10:1$. 
 \setcounter{theorem}{14} 
 \begin{proposition}{\label{3.15}} The projection of the incidence 
 correspondence 
 $$I\subset F'(S)\times 
 VSP(F(S),10)$$ onto the first factor is generically $6:1$. 
 \end{proposition} 
 \proof 
 Consider a general point $g\in F'(S)= P_S\cap K$. Since $g\in K\setminus 
 G$, the 
 Cremona transformation $\varphi$ is defined in $g$ and the fiber containing 
 $g$ is a $5$-space 
 $P_g=P(\wedge^2|g|)$. The space $P_g$ intersects $G$ in a quadric $Q_g$ 
 corresponding to all lines in the 
 $3$-space $P(|g|)$, and $|g|=\ker\varphi(g)$. In the dual space 
 $P_g^{\bot}\subset \Pnd {14}$ is the 
 tangent $8$-space $T_{\varphi(g)}$ to $G^*$ at $\varphi(g)$. Now, the 
 point $p$ is in the image of the 
 projection from the incidence correspondence, i.e. is contained in some 
 $\Gamma\in 
 VSP(F(S),10)$, if $P_g$ is contained in two special tangent hyperplanes, 
 corresponding to two 
 points on $G^*$. These special tangent hyperplanes which contain $P_g$ are 
 parametrized by 
 all $3$-spaces that intersect $P(|g|)$ in a plane. On $G^*$ this is, by 
 \ref{2.5}, precisely what is 
 defined by the minors of a matrix $N_g$ equivalent to $N_{01}$ inside the 
 tangent space 
 $T_{\varphi(g)}=P_g^{\bot}\subset \Pnd {14}$. The hyperplanes which  
 contain both $P_S$ and $P_g$ 
 form $P_S^{\bot}\cap P_g^{\bot}$ in $\Pnd {14}$. This is a $3$-space, since 
 the two $5$-spaces $P_S$ 
 and $P_g$ intersect. Inside $T_{\varphi(g)}$ the minors of $N_g$ define a 
 rational quartic $5$-fold scroll, so the intersection with $L_S=P_S^{\bot}$ 
 is 4 
 points. 
 The secants between the 4 points 
 are secant lines to $L_S\cap G^*$, corresponding to sets $\Gamma$ that contain 
 $g$.\qed\vspace{.2 in}

 \setcounter{subsection}{15} 
 \subsection{\label{3.16}} Finally we turn to the general cubic 4-fold. 
 We want to show that\break $VSP(F(S),10)$ 
 deforms smoothly with $F$, or more precisely that $VSP(F,10)$ for a general 
 cubic 
 4-fold $F$ is a deformation of $VSP(F(S),10)$. But $VSP(F(S),10)/\cong\cF(F')$ 
 is a symplectic manifold, and its space of 
 deformations is smooth of dimension 20 (cf.\cite{BD}). In fact every 
 general member of this 
 deformation is the Fano variety of lines in some 4-fold cubic. Since the 
 variety of cubic 
 fourfolds has dimension 20, it would follow that
 \setcounter{theorem}{16} 
 \begin{theorem}{\label{3.17}} $VSP(F,10)$ for a general cubic fourfold 
 $F$ is isomorphic to the 
 Fano variety of lines in some other cubic fourfold. \end{theorem} 
 \proof 
 Consider the correspondence 
 $$VSP=\{(\Gamma,F)|\Gamma\in VSP(F,10)\}\subset Hilb_o^{10}\Pn 5\times {\bf 
 P}({\rm Sym}^3V)$$ 
 where $\Pn 5={\bf P}(V)$, and $Hilb_o^{10}\Pn 5$ is the component of the 
 Hilbert 
 scheme that contains the smooth subschemes $\Gamma$ that span $\Pn 5$. The 
 fiber of the second 
 projection is $VSP(F,10)$ for a cubic $F\in {\bf P}({\rm Sym}^3V)$. The 
 fiber over a point 
 $\Gamma\in Hilb_o^{10}\Pn 5$ by the first projection is a linear space, the 
 span ${\bf P}(\Gamma)$ of 
 $\rho_3(\Gamma)$ in ${\bf P}({\rm Sym}^3V)$ under the $3-$uple embedding 
 $\rho_3$. For the 
 general $\Gamma$ the span ${\bf P}(\Gamma)$ is of course a $\Pn 9$ , 
 while $Hilb_o^{10}\Pn 5$ has dimension 50, so $VSP$ is reduced of 
 dimension 59. The group $GL(V)$ acts on both factors, and on the 
 incidence $VSP$. Since the 
 general cubic and the general subscheme $\Gamma$ has no nontrivial 
 automorphisms, 
 the quotient by 
 $GL(V)$ is smooth at the general point, and the general fiber of the second 
 projection is 
 the same after taking quotients. Thus for our purposes it suffices to 
 check the second 
 projection in the above incidence. 
 By \ref{3.13} the general fibers of the second projection over the 
 hypersurface in ${\bf 
 P}({\rm Sym}^3V)$ of apolar cubic 4-folds $F(S)$ are smooth 4-dimensional 
 symplectic varieties. Consider a general point $(\Gamma, F(S))$ in one of 
 these fibers.
 \begin{lemma}{\label{3.18}} The projection $VSP\to {\bf P}({\rm Sym}^3V)$ 
 has maximal rank at 
 $(\Gamma, F(S))$.\end{lemma}
 \proof We may assume that $\Gamma$ is smooth, i.e. that 
 $\Gamma=\{l_1,\ldots, l_{10}\}\subset 
 {\bf P}({\rm V})$ and that $f=l_1^3+\ldots +l_{10}^3$, where $F(S)=Z(f)$. 
 The rank of the map 
 onto 
 ${\bf P}({\rm Sym}^3V)$ is then the dimension of the span of 
 $\{l_i^2y_j|1\leq i\leq 10, 0\leq 
 j\leq 5\}$ where $<y_0,\ldots,y_5>=V$. In fact, from the expansion of 
 $(l_i+y_j)^3$, we see that $l_i^2y_j$ defines a tangent direction at the point 
 $l_i^3$, so the above span is the span of the tangent spaces to the $3$-uple 
 embedding $\rho_3({\bf P}(V))$ at the points $l_i^3$. By Terracinis lemma 
 (cf. \cite{Zak}) the tangent 
 space to the $10$th secant variety of $\rho_3({\bf P}(V))$ at $F(S)$ is 
 precisely the span of 
 the tangent spaces at the $10$ points $l_i^3$. If the span is not all of 
 ${\bf P}({\rm 
 Sym}^3V)$, there is a hyperplane section of $\rho_3({\bf P}(V))$ singular 
 at the $10$ points, 
 i.e. a cubic hypersurface in ${\bf P}(V)$ singular in the points 
 $l_1,\ldots, l_{10}$. 
 But it follows from \ref{2.9} 
 that $\{l_1,\ldots, l_{10}\}$ is a general intersection of two rational 
 quartic scrolls 
 in a quadric hypersurface. By genericity it is enough to check in one 
 example that such a 
 subscheme is not the singular locus of any cubic hypersurface. In fact the 
 two matrices 
 \[ 
 \left(\begin{array}{cccc}x_0&x_2&x_4&x_0+x_5\\ 
 x_1&x_3&x_5&x_3+x_4\end{array}\right)\; {\rm 
 and}\; 
 \left(\begin{array}{cccc}x_0&x_1&x_1+x_2&x_1+x_5\\ 
 x_2&x_3&x_4-x_5&x_0+x_3\end{array}\right) 
 \] 
 have a common $2\times 2$ 
 minor. It is easy to check with MACAULAY \cite{MAC} that the two matrices 
 drop rank along 10 points, and 
 that there are no cubic hypersurface singular along these points. 
 Thus the 
 lemma follows. \qed\vspace{.2 in}
 
 The point $(\Gamma, F(S))$ is therefore a regular point for the projection 
 $VSP\to {\bf P}({\rm Sym}^3V)$, i.e. the fibration is smooth at this point. 
 Since the fiber itself 
 is smooth, the general fiber is smooth, and is a deformation of $VSP(F(S),10)$. 
 \qed\vspace{.2 in}
 
 Theorem \ref{3.17} describes a map $\mu$ from an open set of the moduli of 
 cubic 
 $4$-folds into itself. By \ref{1.6} this map is not the identity. 
 Therefore we pose
 \begin{problem} Is $\mu$ dominant, and if so, what is its 
 degree?\end{problem}
 
 \vspace{.2 in}
 
 Authors' addresses:
 \vspace{.2 in}
 
 Atanas Iliev\par 
 Institute of Mathematics, 
 Bulgarian Academy of Sciences,\par 
 Acad. G. Bonchev Str., 8,\par 
 1113 Sofia, Bulgaria.\par 

 e-mail: ailiev@math.bas.bg
 \vspace{.2 in}
 
 Kristian Ranestad\par 
 Matematisk Institutt, UiO,\par 
 P.B. 1053 Blindern,\par 
 N-0316 Oslo, Norway.\par 
 
 e-mail: ranestad@math.uio.no
 \end{document}